\documentclass[12pt]{amsart}
\usepackage[]{hyperref} 
\usepackage{amssymb}
\usepackage[curve]{xypic}
\newcommand{\Omit}[1]{\begin{tiny}#1\end{tiny}}
\renewcommand{\Omit}[1]{}

\newbox\mybox
\def\overtag#1#2#3{\setbox\mybox\hbox{$#1$}\hbox to
  0pt{\vbox to 0pt{\vglue-#3\vglue-\ht\mybox\hbox to \wd\mybox
      {\hss$\scriptstyle#2$\hss}\vss}\hss}\box\mybox}
\def\undertag#1#2#3{\setbox\mybox\hbox{$#1$}\hbox to 0pt{\vbox to
    0pt{\vglue#3\vglue\ht\mybox\hbox to \wd\mybox
      {\hss$\scriptstyle#2$\hss}\vss}\hss}\box\mybox}
\def\lefttag#1#2#3{\hbox to 0pt{\vbox to 0pt{\vss\hbox to
      0pt{\hss$\scriptstyle#2$\hskip#3}\vss}}#1}
\def\righttag#1#2#3{\hbox to 0pt{\vbox to 0pt{\vss\hbox to
      0pt{\hskip#3$\scriptstyle#2$\hss}\vss}}#1}

\def\Dot{\lower.2pc\hbox to 3.5pt{\hss$\bullet$\hss}}
\def\Circ{\lower.2pc\hbox to 3.5pt{\hss$\circ$\hss}}

\def\splicediag#1#2{\xymatrix@R=#1pt@C=#2pt@M=0pt@W=0pt@H=0pt}

\renewcommand\frame[2][3pt]{\hbox{$\vcenter{\hbox{\vrule\vbox
{\hrule\kern#1\hbox{\kern#1$#2$\kern#1}\kern#1\hrule}\vrule}}$}}

\newcommand\lineto{\ar@{-}}
\newcommand\dashto{\ar@{--}}
\newcommand\dotto{\ar@{.}}
\newcommand{\C}{\mathbb C}

\newcommand{\E}{\mathbb E}
\newcommand{\N}{\mathbb N}
\newcommand{\Z}{\mathbb Z}
\newcommand{\Q}{\mathbb Q}

\newtheorem*{DVC}{DCC Volumes Conjecture}

\newtheorem*{theorem*}{Theorem}
\newtheorem{theorem}{Theorem}[section]
\newtheorem{proposition}[theorem]{Proposition}
\newtheorem*{proposition*}{Proposition}
\newtheorem{lemma}[theorem]{Lemma}
\newtheorem{corollary}[theorem]{Corollary}
\theoremstyle{definition}
\newtheorem{example}[theorem]{Example}
\newtheorem*{example*}{Example}
\newtheorem{definition}[theorem]{Definition}
\newtheorem*{definition*}{Definition}
\newtheorem{remark}[theorem]{Remark}
\newtheorem*{remark*}{Remark}

\newtheorem*{conjecture*}{Conjecture}

\begin{document}
\title{The volume of a  surface or orbifold pair}
\author{Jonathan Wahl}
\address{Department of Mathematics\\The University of North
  Carolina\\Chapel Hill, NC 27599-3250} \email{jmwahl@email.unc.edu}
  \begin{abstract}A \emph{surface pair} $(X,C)$ is a germ of a complex normal surface singularity $(X,0)$ and a sum $C=\sum c_iC_i$ of curves on $X$, with $c_i\in [0,1]\cap\Q$.  An \emph{orbifold pair} has $c_i=1/n_i$, as intersecting with a small sphere gives a $3$-dimensional orbifold $(\Sigma, \gamma_i,n_i)$.  There are natural notions of morphism and log cover of surface pairs.  We introduce a volume $Vol(X,C)\in \Q_{\geq 0}$, computable from any log resolution, analogous to that in our 1990 JAMS paper when $C=0$.  Denoting $\bar{C}=\sum (1-c_i)C_i$, one has $(X,C)$ log canonical iff $Vol(X,\bar{C})=0$.  The main theorem (5.4-5.5) is that $Vol(X,C)$ is ``characteristic'': if $f:(X',C')\rightarrow (X,C)$ is a morphism of degree $d$, then $Vol(X',C')\geq d\cdot Vol(X,C)$, with equality if $f$ is a log cover.  We prove (6.7) that $Vol(X,\sum (1/n_i)C_i)=0$ iff the associated orbifold has finite or solvable orbifold fundamental group, and these are classified.  In (8.2) is proved a key case of the DCC Volumes Conjecture: The set $\{Vol(X,\sum (1/n_i)C_i)|X\  \text{RDP}\}$ satisfies the DCC, with minimum non-$0$ volume 1/3528. 
  \end{abstract}
\maketitle

\section{{Introduction}}

Milnor and Thurston \cite{MT} defined a \emph{characteristic number} to be a real-valued function $\lambda$ on a class of compact $3$-manifolds with the property that it multiplies by degree in an unramified covering.   As an example, for the class of compact hyperbolic $3$-manifolds, compatibility of the metrics in a finite covering means \emph{volume} is a characteristic number.  In \cite{MT} they defined such an invariant $\lambda$ for all $3$-manifolds with a further important submultiplicative property: if $f:M'\rightarrow M$ has degree $d$, then $\lambda (M')\geq d \lambda (M)$, with equality if $f$ is a covering. However, this invariant is extremely difficult both to define and to compute.

In our paper \cite{w}, we defined a characteristic number for the class of links of complex normal surface singularities (NSS's).  Intersecting the germ $(X,0)$ of a NSS with a small sphere yields its link $\Sigma$, a compact oriented $3$-manifold.  For a good resolution $(\tilde{X},E)\rightarrow (X,0)$, $\Sigma$ is the boundary of a tubular neighborhood of $E$ on $\tilde{X}$, and its topology is completely determined by the weighted dual graph $\Gamma$ of $E$.  

 Using the Zariski decomposition $K_{\tilde{X}}+E=P+N$ of $\Q$-divisors on any good resolution, we showed $-P\cdot P$ is a characteristic number.   It is now renamed as $Vol(X)$. 
\begin{theorem}\cite{w}  For a NSS $(X,0)$, there is an invariant $Vol(X)$, computable from any resolution graph $\Gamma$, for which:
\begin{enumerate}
\item $Vol(X)\in \Q_{\geq 0}$ is a characteristic number for the link $\Sigma$.
\item If $f:(X',0)\rightarrow (X,0)$ is finite of degree $d$, one has \begin{center}$Vol(X')\geq d \cdot Vol(X),$
\end{center} with equality if $f$ is unramified off $0$.
\item $Vol(X)=0$ if and only if $(X,0)$ is log canonical if and only if the local fundamental group $\pi_1(\Sigma)$ is finite or solvable. 
\item If $(X,0)$ is quasi-homogeneous and not log canonical, there is a natural metric on the Seifert link $\Sigma$ whose volume is $4\pi^2Vol(X)$.
\end{enumerate}
\end{theorem}
Note that ramified maps could have equality in the conclusion of $(2)$, for instance considering any map of degree at least 2 of $(\C^2,0)$ into itself.

 In \cite{B},  Boucksom, de Fernex, and Favre define a volume for higher dimensional isolated normal singularities, proving property (2); the definition was much more complicated, requiring deep work of Shokurov.

In analogy with results about the set of volumes of hyperbolic manifolds, the following theorem was conjectured in \cite{w} and proved by F. Ganter.

\begin{theorem}\cite{G} The set of volumes $\{Vol (X)|(X,0)  \text {Gorenstein}\}$ satisfies the DCC, with smallest non-$0$ value $1/42$, obtained precisely for the Brieskorn singularity $V(2,3,7)$ and its equisingular deformation.
\end{theorem}

The goal of the present work is to extend the results of \cite{w} to \emph{surface pairs} $(X,C=\sum_{i=1}^r c_iC _i)$, where $(X,0)$ is a NSS and the ``curve'' $C$ consists of irreducible Weil divisors $C_i$ with $c_i\in [0,1]\cap \Q$.  Write $C'\leq C$ if $C-C'$ is an effective divisor.  We say $C$ is \emph{reduced} if all the non-zero $c_i$'s are equal to $1$.

An important special case are the \emph{orbifold pairs} $(X,\sum_{i=1}^r (1/n_i)C _i)$, with $n_i\in \N$ (regrettably, called ``singular pairs'' in \cite{nw} and written $(X,\sum_{i=1}^r n_iC_i$)). As  remarked in \cite{nw}, intersecting $(X,0)$ and the $C_i$ with a small sphere gives a $3$-dimensional  \emph{orbifold} $(\Sigma, \gamma_1,\cdots,\gamma_r, n_1,\cdots,n_r)$, where the $\gamma_i$ are knots in $\Sigma$.  

Surface pairs are studied via a \emph{log resolution} $\pi:(\tilde{X},\tilde{C}\cup E)\rightarrow (X,C)$, in which $\tilde{C}\cup E=\pi^{-1}(\cup C_i)$ has normal crossings.  Calculations can be made from an associated decorated graph $\Gamma^*$, as in \cite{nw}, showing the location of the proper transforms $\tilde{C_i}$ of the $C_i$ and the associated weights $c_i$ (or $n_i$).   Consider the Zariski decomposition on $\tilde{X}$: $$K_{\tilde{X}}+E+\sum c_i\tilde{C}_i=P_{\tilde{C}}+N_{\tilde{C}}.$$  In Section 3 are established basic properties:

\begin{proposition*}[\ref{a}]
\begin{enumerate}
\item $P_{\tilde{C}}=0$ iff $(X,C)$ is log canonical.
\item $-P_{\tilde{C}}\cdot P_{\tilde{C}}:= -P_C\cdot P_C\in \Q_{\geq 0}$ is independent of the log resolution.
\item $C'\leq C$ implies $-P_{C'}\cdot P_{C'}\leq -P_C\cdot P_C$.
\item $C=0$ implies $-P_C\cdot P_C=Vol(X)$.
\item For $X$ log terminal and $C$ reduced, the log canonical threshold of $C$ is the maximum $\epsilon$ with $-P_{\epsilon C}\cdot P_{\epsilon C}=0$
\end{enumerate}
\end{proposition*}

An \emph{end curve} of $E$ is an irreducible component of $E$ which intersects only one other such component.  For calculations, it is easiest to use the \emph{minimal orbifold log resolution} of a pair, with the property that every $\tilde{C_i}$ intersects an end curve of $E$, and each end intersects at most one $\tilde{C_i}$.  In this case, if $\Gamma^*$ is star-shaped, one can easily write a formula for $-P_C\cdot P_C$ in terms of weights and lengths of chains (Proposition $4.2$).  It is important to be able to compute this invariant when all $c_i\geq 1/2$, especially for calculations of orbifold volumes (Section 6). 
\begin{theorem*}[\ref{t}] Consider a pair $(X,\sum c_iC_i)$ for which all non-$0$ $c_i$ are $\geq 1/2.$  If $\Gamma^*$ is not star-shaped, then $-P_C\cdot P_C$ can be computed directly from the minimal orbifold log resolution, without performing a Zariski decomposition.
\end{theorem*}

In Section 5, we define maps and coverings of pairs.  In \cite{nw}, an orbifold covering of $(\Sigma, \gamma_i,n_i)$ was defined as a map of an orbifold to $\Sigma$ which is branched only over the $\gamma_i$, and with branching order $\leq n_i$.  Here we adapt that notion for maps of surface pairs.  For  a pair $(X,C)$ and map $f:X'\rightarrow X$, one can define (Section 5) a pair $(X',f^{*}(C))$ which encodes the $c_i$ and the branching orders of the components of the inverse image of $|C|$.
\begin{definition*}[\ref{b}]  A \emph{map of surface pairs} $f:(X',C')\rightarrow (X,C)$ is a finite map $f:(X',0)\rightarrow (X,0)$ for which $f^*(C)=C'$.  A \emph{log cover} is a map of surface pairs which is unramified off the support of $C'$.
\end{definition*}
For any curve $C=\sum c_iC_i$, the \emph{opposite curve} is $\bar{C}:=\sum (1-c_i)C_i$.   The crucial definition is:
\begin{definition*}[\ref{c}]  For $(X,C=\sum c_iC_i)$ a surface pair, the \emph{volume of the pair} is $$Vol(X,C):=-P_{\bar{C}}\cdot P_{\bar{C}}.$$
\end{definition*}

The major result of this paper, proved primarily in (\ref{g}) and (\ref{h}), is
\begin{theorem*}[\ref{g}),(\ref{h}]  Let $(X,C=\sum c_iC_i)$ be a surface pair, with decorated resolution graph $\Gamma^*$ and all $c_i\neq 0$.
\begin{enumerate}
\item  $Vol(X,C)\in \Q_{\geq 0}$ is computable from $\Gamma^*$.
\item If $f:(X',C')\rightarrow (X,C)$ is a map of pairs of degree $d$, then 
\begin{center} $Vol(X',C')\geq d\cdot Vol(X,C)$,
\end{center} with equality if $f$ is a log cover.
\item If $C$ is reduced (all $c_i=1$), then $Vol(X,C)=Vol(X)$.
\item $C'\leq C$ implies $Vol(X)\leq Vol(X,C)\leq Vol(X,C')$.
\end{enumerate}
\end{theorem*}
Restricting to the class of orbifolds arising from orbifold pairs $(X,C)$, Theorems $5.4$ and $5.5$ imply that $Vol(X,C)$ is a characteristic number in the sense of Milnor-Thurston: multiplicative by degree for orbifold covers, submultiplicative for orbifold maps.
\begin{theorem*}[\ref{k}] Suppose $f:(X',\sum (1/n'_j)C'_j)\rightarrow (X,\sum (1/n_i)C_i)$ is a map of degree $d$ of orbifold pairs. Then $$Vol (X',C')\geq d\cdot Vol(X,C),$$ with equality if $f$ is a log cover.
   \end{theorem*}
Recall that orbifold covers of an orbifold $(\Sigma, \gamma_i,n_i)$ are controlled by an \emph{orbifold fundamental group} $\pi_1^{orb}(\Sigma)$ (6.3).  Using Theorem $6.2$ and the notion of \emph{UALC (universal abelian log cover)} from \cite{nw}, we prove a result analogous to the classification of log canonical singularities:
\begin{theorem*} [\ref{l}] The volume of an orbifold pair is $0$ (i.e., it  is log canonical) if and only if the corresponding orbifold fundamental group is finite or solvable.  
\end{theorem*}
Via computations of Section 4 and consideration of star-shaped graphs $\Gamma^*$ with $3$ branches, one has:
\begin{theorem*} [\ref{m}] The resolution graphs of all volume $0$ orbifold pairs are classified.
\end{theorem*}
An orbifold pair $(X,\sum (1/n_i)C_i)$ has volume $0$ iff $(X,\sum(1-1/n_i)C_i)$ is log-canonical.  Since J. Koll\'{a}r (\cite[(3.3)]{k}) has classified log canonical surface pairs with all non-zero $c_i\geq 1/2$ , the last result could be deduced from his list.  In fact, the approach of this paper can yield his more general result as well ((7.5)-(7.8)).

Simple examples indicate that the set of volumes of Gorenstein orbifold pairs does not satisfy the DCC.  But the following seems to be a reasonable 
\begin{DVC}
The set of volumes of orbifold pairs $(X,C)$ on rational double points satisfies the DCC, and has minimum non-$0$ value $1/3528=(2\cdot 42\cdot 42)^{-1}.$
\end{DVC}
It makes sense to consider first the non-trivial case of star-shaped $\Gamma^*$.
\begin{theorem*} [\ref{n}] Consider the set of volumes of Gorenstein orbifold pairs $(X,C)$ for which $(X,0)$ is a rational double point and the minimal orbifold log resolution is star-shaped.  Then this set satisfies the DCC, with  minimum non-$0$ volume $1/3528$.
\end{theorem*}

At present, we see no metric interpretation of $Vol(X,C)$ for an orbifold pair when $(X,0)$ and the $C_i$ are quasi-homogeneous.

In a separate work, writing $C=\sum c_iC_i$, we shall discuss the piecewise-quadratic behavior of the function $-P_C\cdot P_C$ as a function of $(c_1,\cdots,c_r)$, with $\{(c_1,\cdots,c_r)|\  0\leq c_i\leq 1\}$ divided into relevant regions by linear inequalities.  This is particularly relevant for $cC$, where $C$ is a reduced plane curve germ.

Wim Veys has pointed out an alternative way to view the invariant $-P_C\cdot P_C$, which was here defined using Zariski decomposition $P+N$ of $K_{\tilde{X}}+E+\sum c_i\tilde{C}_i$ on a log resolution $(\tilde{X},E)$ of $(X,C)$.  Instead, pass to the (unique) log canonical model $(X_{lc},E_{lc})$ of $(X,C)$, obtained for instance by blowing down $\tilde{X}$ along $N$ and any other $E_i$ dotting to $0$ with $P$.
Here the corresponding pair $(X_{lc},E_{lc}+\sum c_i\tilde{C_i})$ has only log canonical singularities and $K_{X_{lc}}$ +$E_{lc}+\sum c_i\tilde{C}_i$ is relatively ample.  Veys shows that $$-P_C\cdot P_C=-(K_{X_{lc}} +E_{lc}+\sum c_i\tilde{C}_i)^2.$$  One could proceed similarly using the (unique) minimal dlt resolution $(X_{dlt},E_{dlt})$ of $(X,C)$.

We are very grateful to Enrique Artal Bartolo for several important comments which greatly improved this paper, and to the organizers (including NSF) of the 115AM Jaca Conference in June, 2023.

\section{{Review of Zariski decomposition and \cite{w}}}

Let $(\tilde{X},E)\rightarrow (X,0)$ be a good resolution of a normal surface singularity.  The adjoint homomorphism $$\text{Pic}\  \tilde{X}\rightarrow \bigoplus \Q\cdot E_i=: \mathbb E_\Q,$$ sends a line bundle $\mathcal L$ to the exceptional $\Q$-divisor $\Sigma a_iE_i$ satisfying $$\mathcal L\cdot E_j=(\Sigma a_iE_i)\cdot E_j, \ \text{all} \ j.$$  Denote the image of $\mathcal L$ by $L$, while $K$ or $K_{\tilde X}$ denotes both the canonical line bundle and its image; similarly, use the same notation for a curve $\tilde{C}$ on $\tilde{X}$ and its image under the map. 
 By abuse of notation, we shall write equality ($=$) when we really mean numerical equivalence ($\equiv$) of $\Q$-divisors.

\begin{proposition} (Sakai  \cite[p. 408]{s}) Let $\L$ be an exceptional $\Q$-divisor.  Then there exists a unique Zariski decomposition $L=P+N$ in $ \mathbb E_{\mathbb Q}$, where
\begin{enumerate}
\item $P$ is nef, i.e., $P\cdot E_i\geq 0$, all $i$.
\item $N$ is effective, i.e., $N=\Sigma a_iE_i$, with all $a_i\geq 0$.
\item $P\cdot N=0$, i.e., $P\cdot E_i=0$ for all $E_i\subset \text{Supp} \ N$.
\end{enumerate} 

\end{proposition}

\begin{remark} The usual inductive approach to constructing $N$ from $L$ starts with letting $\mathcal A_1=\{j| \ L\cdot E_j<0\}$, and defining $N_1=\Sigma b_jE_j$ ($j\in \mathcal A_1$) so that $L\cdot E_j=N_1\cdot E_j$, for all $j\in \mathcal A_1$.  If $P_1=L-N_1$ is nef, we are done; if not, $\mathcal A_2$ and $N_2$ are needed.  This would occur if a neighbor $E_k$ of the support of $N_1$ dots to $0$ with $L$, which curve would have to be included in $N_2$.  So, an ``improved'' $N_1$ can be defined as follows:   If $L$ is not nef, consider the union of all $E_i$ such that $L\cdot E_i\leq 0$, and pick out the connected components which contain a curve with $L\cdot E_j<0$.  Let  $\mathcal A_1$ correspond to all the $E_i$ in such components.  Choose $N_1$ as above using this $\mathcal A_1$, so that $P\cdot E_k=N_1\cdot E_k$ for all $k\in \mathcal A_1$.  Often, $P_1=L-N_1$ will be nef. If not, continue as above to find $N_2$.

\end{remark}

\begin{remark}  The Zariski decomposition $P+N$ of $L$ can be constructed directly from the knowledge of Supp $N$.  For, $N$ is the linear combination of the curves in its support satisfying $L\cdot E_i=N\cdot E_i$ for all such $E_i$. 
\end{remark}
To prove Theorems $5.4$ and $5.5$ below, we
use some results from \cite[(1.4),(1.6)]{w} :
\begin{enumerate} 
\item Zariski decomposition is preserved under pullback of a generically finite and proper surjective map $f:(\tilde{Y},F)\rightarrow (\tilde{X},E)$, since $f^*:\mathbb E_\Q\rightarrow \mathbb F_Q$ preserves both effective and nef divisors, and multiplies intersection numbers by deg($f$).
\item Suppose $\mathcal L$ is a line bundle on $\tilde{X}$, with $L=P+N$ the associated Zariski decomposition of $\Q$-divisors.  Then for $n>>0$, $$\text{dim}\ H^0(\tilde{X}-E, \mathcal L^n)/H^0(\tilde{X},\mathcal L^n)=n^2/2(-P\cdot P)+O(n).$$ 
\item Suppose $\mathcal L':=\mathcal L(-D)$, where $D$ is an effective exceptional divisor.  If $L'=P'+N'$, then $$-P\cdot P\leq -P'\cdot P'.$$
\end{enumerate}
The third statement follows from the second because $H^0(\tilde{X},\mathcal L'^n)\subset H^0(\tilde{X},\mathcal L^n).$

We need a slight generalization of the third statement.
\begin{proposition}  Let $L=P+N$ and $L'=P'+N'$ be Zariski decompositions of two $\Q$-divisors.  Assume there is an effective $\Q$-divisor $D$ with $L'=L-D.$  Then $$-P\cdot P\leq -P'\cdot P'.$$
\begin{proof}   Multiplying a $\Q$-divisor by an integer $k$ multiplies its $-P\cdot P$ by $k^2$.  So we can multiply $L$ and $D$ by $k$ and assume the divisors $L$ and $L'$ are integral and $D$ is integral and effective.   Define $\mathcal L=\mathcal O(L)$ and $\mathcal L'=\mathcal L(-D)$ and apply the third statement above.
\end{proof}
\end{proposition}
In \cite{w} we considered the Zariski decomposition $K_{\tilde {X}}+E=P+N$ on a good resolution $\tilde{X}$ of $(X,0)$, and we showed  $-P\cdot P=: Vol(X)$ is a non-negative rational number independent of the resolution.  On the minimal good resolution of $(X,0)$, it is straightforward to write down $N$, thus enabling the calculation of $Vol(X)$; this allowed recovery of Kawamata's classification in \cite{K} of the
log canonical singularities, where $P=0$ (i.e., $K_{\tilde{X}}+E$ is numerically effective). 
The key properties of $Vol(X)$ are listed in Theorem 1.1 in the Introduction.

\section{{Surface pairs and their invariant $-P\cdot P$ }}
     
     We consider a surface pair $(X,C)$, where  $(X,0)$ is a NSS and $C=\sum_{i=1}^r c_iC_i$ is a finite sum of irreducible Weil divisors, with $c_i\in  [0,1]\cap\Q$.  If $C'=\sum_{i=1}^r c'_iC_i$, write $C'\leq C$ if  $c'_i\leq c_i$ for all $i$.  $C=0$ means all $c_i=0$; $C$ \emph{reduced} means all non-$0$ $c_i$ equal $1$.   A \emph{log resolution} of  $(X,C)$ is a good resolution $\pi:(\tilde{X},E)\rightarrow (X,0)$ for which $E\cup \pi^{-1}(\cup C_i)$ has strong normal crossings.  Denote the proper (sometimes called birational, or strict) transform of $C_i$ on
      $\tilde{X}$ by $\tilde{C}_i$. As each
      $\tilde{C}_i$ dots non-negatively with every $E_j$, and positively with at least one such, we conclude that $\tilde{C}_i$ is numerically equivalent to a $\Q$-divisor all of whose coefficients are negative rational numbers.

     Writing $$K_{\tilde{X}}+\sum c_i\tilde{C}_i\equiv \sum a_i E_i,$$  the $a_i$ are 
the \emph{discrepancies} of the pair $(X,C)$.  Often one assumes $K$ and the $C_i$ are $\Q$-Cartier, but this is not necessary.

      $(X,C)$ is called \emph{log canonical} if all $a_i\geq -1$, \emph{log terminal} if all $a_i>-1.$ 
      
      As in \cite{w}, it is natural to consider instead the $\Q$-divisor $$K_{\tilde{X}}+E+\sum c_i\tilde{C}_i$$  and its Zariski decomposition  
      $$K_{\tilde{X}}+E+\sum c_i\tilde{C}_i=P_{\tilde{C}}+N_{\tilde{C}}.$$
      Then log canonicality is equivalent to $P_{\tilde{C}}=0$, and log terminality to the additional requirement that $N_{\tilde{C}}$ is strictly effective (i.e., all coefficients are positive).
\    
     \begin{proposition}  \label{a}The non-negative rational number $-P_{\tilde{C}} \cdot P_{\tilde{C}}$ is independent of the log resolution of the pair $(X,C)$.
     \begin{proof} (Cf.\cite[(2.7)]{w}).  Let $\pi:(X',F)\rightarrow (\tilde{X},E)$ be the blow-up along a point $q$ of $E$, adding an exceptional $-1$-curve $F_1$ to $F$.  Then one has $$K_{X'}+F=\pi^*(K_{\tilde{X}}+E)+\delta F_1,$$ where $\delta =1$ if $q$ is a smooth point of $E$, and $\delta=0$ if $q$ is a double point.  Further, if $C'_i$ is the proper transform of $\tilde{C}_i$, then $$C'_i=\pi^*(\tilde{C}_i)-\epsilon_i F_1,$$ where $\epsilon_i=1$ if $q\in \tilde{C}_i$ and is $0$ otherwise. 
      Putting these together gives $$K_{X'}+F+\sum c_iC'_i=\pi^*(K_{\tilde{X}}+E+\sum c_i\tilde{C}_i)\ +\ (\delta-\sum\epsilon_i c_i)F_1.$$  Note that the coefficient of $F_1$ on the right is non-negative in all cases (e.g., equal $1-c_i\geq 0$ if $q\in \tilde{C}_i$, equal to $0$ if $q$ is a double point).  Since $\pi^*(P_{\tilde{C}})\cdot F_1=0$ and Zariski decomposition is preserved under pull-back, one has on $X'$ a new Zariski decomposition
       $$K_{X'}+F+\sum c_iC'_i=\pi^*(P_{\tilde{C}}) +(\pi^*(N_{\tilde{C}})+(\delta-\sum\epsilon_i c_i)F_1).$$ 

    It follows that that $\pi^*(P_{\tilde{C}})=P_{C'}$, whence the assertion.
     \end{proof}
     \end{proposition}
     
 \begin{definition}  For a surface pair $(X,C)$, let  \emph{$-P_C \cdot P_C$} be the quantity $-P_{\tilde{C}} \cdot P_{\tilde{C}}$ computed from any log resolution.
 \end{definition}
  
  It will be convenient to consider only log resolutions of $(X,C)$ for which every $\tilde{C_i}$ intersects an end curve of $E$, with at most one $\tilde{C_i}$ per end-curve.  As in \cite[(6.4)]{nw}, there exists a unique minimal such log resolution (except when $E$ consists of one curve and there are exactly two $\tilde{C_i}$, intersecting it at two points), called there the \emph{minimal orbifold resolution}. We can restrict our attention to this resolution.   In this case, the line bundle defined by $\tilde{C}_i$, which intersects only $E_i$, is numerically equivalent to minus the effective $\Q$-divisor  $e_i$ satisfying $e_i\cdot E_j=-\delta_{ij}$.
  \begin{proposition}   
  $C'\leq C$ implies $-P_{C'} \cdot P_{C'} \leq -P_{C} \cdot P_{C}$.
  \begin{proof}  Adding a multiple of $\tilde{C_i}$ is equivalent to subtracting an effective divisor, so the result follows from Proposition 2.4.
  \end{proof}
   \end{proposition}
  The log canonical pairs with $C$ reduced ($c_i=0$ or $1$) were classified completely by Kawamata \cite{K},  as previously mentioned.  J. Koll\'{a}r has listed (\cite[(3.3)]{k})the examples with all $c_i\geq 1/2$, and we indicate how to rederive this result in $((7.5)-(7.8))$ below.   Calculations are much harder for smaller non-$0$ values of $c_i$. 
  
  Because each $\tilde{C}_i$ is equivalent to a strictly negative divisor, one concludes
  
  \begin{lemma}
  \begin{enumerate}
  \item If $(X,C)$ is log canonical, and $C' <C$,  then $(X,C')$ is log terminal.
  \item If $(X,0)$ is log terminal, $C$ arbitrary, then for sufficiently small $\epsilon>0$, $(X,\epsilon C)$ is log terminal.
  \end{enumerate}
  \begin{proof} $(X,C)$ log canonical means $K_{\tilde{X}}+E+\tilde{C}$ is equivalent to an effective divisor, and subtracting the strictly negative divisor $\tilde{C}-\tilde{C'}$ gives a strictly positive divisor, hence $(X,C')$ is log terminal. If $D_1$ and $D_2$ are strictly positive divisors, then for all small $\epsilon$, $D_1-\epsilon D_2$ is strictly positive; this gives the second statement.
  \end{proof}
  \end{lemma}
  \begin{definition} If $(X,0)$ is log terminal, C a reduced curve, then the \emph{log canonical threshold (or LCT) of $(X,C)$} is the largest $\epsilon\leq 1$ for which $(X,\epsilon C)$ is log canonical.
  
  \end{definition}
  
  \section{{Examples with simplest graphs}}

  \begin{lemma} Suppose $(X,C)$ has a log resolution $\tilde{X}$ for which $E$ is a chain of smooth rational curves, and $\tilde{C}_i$ ($i=1, 2$) intersect the ends transversally:
  $$\xymatrix@R=4pt@C=24pt@M=0pt@W=0pt@H=0pt{\\
\tilde{C_1}\leftarrow&\lineto[l]\undertag{\bullet}{}{8pt}\lineto[r]&\undertag{\bullet}{}{8pt}\dashto[r]&\dashto[r]&\undertag{\bullet}{}{4pt}\lineto[r]&\rightarrow\tilde{C_2}\\&~}
$$    
    Then $(X, c_1C_1+c_2C_2)$ is log terminal, except that when $c_1=c_2=1$ it is only log canonical.
  \begin{proof}  If $D$ is an exceptional $\Q$-divisor satisfying $D\cdot E_i\leq 0$ for all $i$, then $D$ is either $0$ or is strictly positive. Dotting $K_{\tilde{X}}+E+\sum_{i=1}^2c_i\tilde{C}_i$ with an interior curve is $0$; dotting with an end gives $c_i-1$ if the chain is not a single curve, or $c_1+c_2-2$  in case $E$ is one curve.
  \end{proof}
  \end{lemma}
    From now one we shall pretty much exclude this simple case from consideration. 
 
    To find Zariski decomposition on the minimal orbifold resolution $(\tilde{X},\tilde{C}\cup E)\rightarrow (X,C)$, it is clear from \cite[(2.3)]{w} one must consider maximal end strings of rational curves $E_1,\cdots,E_s$ contained in the graph $\Gamma$.  That is,  $E_i\cdot E_{i+1}=1,\  1\leq i <s$; $E_2,\cdots, E_s$ have no other intersections in $E$ (so $E_s$ is an end curve); and $E_1$ intersects exactly one other curve, either of positive genus or valence at least $3$.   Writing $a_j=-E_j \cdot E_j$, one has a continued fraction $m/q=a_1-1/a_2-\cdots -1/a_s$; by minimality, all $a_j\geq 2$, unless the string consists of one $-1$ curve  (when $m=1$).  The expansion going in the other direction, starting from $a_s$, corresponds to $m/q'$, where $qq' \equiv 1$ mod $m$. To an $m/q$ string we associate effective $\Q$-cycles:
\begin{itemize}
\item  $D=\sum_{j=1}^s  b_j E_j$, with $D\cdot E_s=-1$ and $D\cdot E_j=0, j<s$; then  $b_1=1/m$ and $b_s=q'/m$.
\item  $D'=\sum_{j=1}^s b'_jE_j$, with $D'\cdot E_1=-1$ and  $D'\cdot E_j=0, j>1$; then $b'_1=q/m$ and $b'_s=1/m.$
\end{itemize}

Since $K+E$ dots to $-1$ with the end-curve $E_s$ and to $0$ with other curves in the string, as indicated in Remark $2.2$ the $N$ in its Zariski decomposition must include the cycle $D$.

          Consider a star-shaped resolution graph, consisting of: a smooth central curve $F$ of genus $g$ and self-intersection $-d$, and $t$ strings $S_i$ with invariants $m_i/q_i$, emanating from $F$.
          $$\xymatrix@R=4pt@C=24pt@M=0pt@W=0pt@H=0pt{\\
\lefttag{\bullet}{m_2/q_2}{8pt}\dashto[ddrr]&
&\hbox to 0pt{\hss\lower 4pt\hbox{.}.\,\raise3pt\hbox{.}\hss}
&\hbox to 0pt{\hss\raise15pt
\hbox{.}\,\,\raise15.7pt\hbox{.}\,\,\raise15pt\hbox{.}\hss}
&\hbox to 0pt{\hss\raise 3pt\hbox{.}\,.\lower4pt\hbox{.}\hss}
&&\righttag{\bullet}{m_{t-1}/q_{t-1}}{8pt}\dashto[ddll]\\
\\
&&\bullet\lineto[dr]&&\bullet\lineto[dl]\\
\lefttag{\bullet}{m_1/q_1}{8pt}\dashto[rr]&&
\bullet\lineto[r]&\overtag{\circ}{-d}{8pt}\undertag{}{}{6pt}\lineto[r]&\bullet
\dashto[rr]&&\righttag{\bullet}{m_{t}/q_{t}}{8pt}\\
&&&F&~\\}
$$
        
   Consider the two topological  invariants introduced (except for a change of sign) by Neumann \cite[p. 250]{n} :  $$\chi:=2g-2+t-\sum_{i=1}^t(1/m_i)$$
          $$e:=d-\sum_{i=1}^t(q_i/m_i).$$
         A $-1$ string (with $m_i=q_i=1$) collapsed to a point leaves $\chi$ and $e$ unchanged.  Recall $e>0$; $\chi<0$ exactly for quotient singularities; $\chi=0$ only for the familiar log canonical star-shaped graphs; and in the Gorenstein case there is an integer $t$ with $\chi=te$  (e.g., \cite[(4.3)]{W}). It is easy to see that if $\chi>0$, then $\chi\geq 1/42$, which occurs for $g=0, t=3, \{m_1,m_2,m_3\}=\{2,3,7\}$.
         
          Noting that $F\cdot D'_i=q_i/m_i$, one has that $F+\sum D'_i$ dots with $F$ to $-e$ and to $0$ with all other curves.   Thus $f$, the effective $\Q$-cycle dotting to $-1$ with $F$ and $0$ with other curves,  satisfies  
 \begin{equation*}\tag{3.5.1}          F+\sum_{i=1}^tD'_i\equiv e f.\end{equation*}
           One concludes that $f\cdot f = -1/e$ (which is why $e>0$).
          
          Suppose now that $(X,\sum_{i=1}^t c_iC_i)$ has star-shaped minimal orbifold log resolution $(X,E)$ as above, where $\tilde{C_i}$ intersects the end of the $i^{th}$ string.  One has a decorated graph $\Gamma^*$, where $\tilde{C_i}$ is represented by an arrow emerging from the end of the string, with a $c_i$ marked over it (e.g.,  (4.5) below).  Assume at least one $c_i\neq 0$, so $t\geq 1$; $m_i=1$ implies $c_i\neq 0$; and $t\geq 3$ if $g=0$. 
      \begin{proposition}  With $(X,\sum_{i=1}^t c_iC_i)$ as above, 
      define 
       $$\chi_C:=\chi+\sum_{i=1}^t c_i/m_i.$$
      Then $(X,C)$ is 
      \begin{enumerate}
      \item  log terminal iff $\chi_C<0$;
      \item  log canonical iff $\chi_C\leq 0$;
      \item not log canonical iff $\chi_C>0$, in which case  $$-P_C\cdot P_C=(\chi_C)^2/e.$$
      \end{enumerate}
      \begin{proof}  If $D_j$ is the cycle above associated to the string $S_j$, then $(1-c_j)D_j$ is effective and dots to $0$ with $K+E+\sum c_i\tilde{C}_i$; further,  $D_j\cdot F=1/m_j.$  With $f$ as above, one checks that $$K+E+\sum c_i\tilde{C}_i=(-\chi_C) f+\sum (1-c_i)D_i,$$ as seen by dotting each side with an end of $S_i$, an interior curve of $S_i$,  and with $F$.  $f$ is an effective cycle, and is the only cycle on the right whose support contains $F$.
      
      The assertions about $\chi_C<0$ and $\chi_C\leq 0$ now follow, because the right-hand side is effective (or $0$).   
      
      If $\chi_C>0$, then $(-\chi_C)f$ is a positive cycle, so equals $P_{\tilde{C}}$, the positive part of the Zariski decomposition of  $K+E+\sum c_i\tilde{C}_i$.  As $f\cdot f=-1/e$, one has $-P_C\cdot P_C=(\chi_C)^2/e$.

      \end{proof}
      \end{proposition}
      
     For log canonicality of $(X,C)$  $(C\neq 0)$ in star-shaped cases,  one must have $(X,0)$ log terminal, hence a quotient singularity (or smooth).  Then $-P_C\cdot P_C$, as a function of the $c_i$'s, is $0$ in the region $\sum c_i/m_i \leq -\chi$ (and $0\leq c_i\leq 1$, all $i$), and given by  the quadratic polynomial $(\chi +\sum c_i/m_i)^2$ outside that region.  Note log canonical  examples can have many strings of $-1$ curves, but at most $3$ whose determinant (i.e.,  $m_i$) is greater than $1$.  One easily deduces the following:
      \begin{corollary} Assume $(X,C)$ as above has all $c_i$ equal to $0$ or $1$.  If $(X,C)$ is log canonical, then the graph is
$$
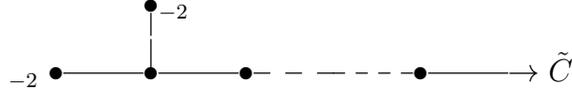
\xymatrix@R=8pt@C=30pt@M=0pt@W=0pt@H=0pt{
&&&\righttag{\bullet}{-2}{6pt}&&&&&&&\\
&&&\lineto[u]&&&&&&\\
&&\lefttag{\bullet}{-2}{4pt}\lineto[r]&\undertag{\bullet}{}{4pt}\lineto[r]\lineto[u]&\undertag{\bullet}{}{4pt}\dashto[r]&\dashto[r]&\undertag{\bullet}{}{4pt}\lineto[r]&\rightarrow\tilde{C}\\
&&&&&\\
}
$$       
      \end{corollary}
   
\begin{corollary} Assume $(X,0)$ is log terminal.  Then excluding the example of Corollary 4.3, the log canonical threshold of the reduced curve $C=\sum_{i=1}^{t'} C_i$ (for some $t' \leq t$) is $$[2-\sum_{i=1}^t (1-1/m_i)]/\sum_{i=1}^{t'}(1/m_i).$$
\end{corollary}

The complete list of examples for which the log canonical threshold is at least $1/2$ can be determined from Proposition $7.7$ below.
 \begin{example} Consider the pair $(\C^2,\text{c}C)$, where $(C,0)$ has an ordinary cusp singularity.  Recalling that an unmarked vertex is a rational $-2$ curve, the decorated graph for the minimal orbifold log resolution is 
 $$
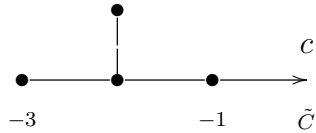
\xymatrix@R=8pt@C=30pt@M=0pt@W=0pt@H=0pt{
&&&&\righttag{\bullet}{}{6pt}&&&&\\
&&&&\lineto[u]&&c&&\\
&&&\undertag{\bullet}{-3}{4pt}\lineto[r]&\overtag{\bullet}{}{6pt}\lineto[r]\lineto[u]&\undertag{\bullet}{-1}{4pt}\lineto[r]&\undertag{\ar}{\tilde{C}}{8pt}\\
&&&&&&\\
&&&&&\\
}
$$  
Since $\chi=-5/6$, one has log canonicality only up to $(5/6)C$.   For $5/6\leq c\leq 1$, $-P_C\cdot P_C$ equals $6(c-5/6)^2$.
 \end{example} 
 
 In the non-star-shaped case with $C\neq0$, the minimal orbifold log resolution $(\tilde{X},\tilde{C}\cup E)$ of $(X,C)$ has maximal end strings, hence associated cycles $D_j$ and $D'_j$.  The effective divisor $N_1=\sum (1-c_j)D_j$ is orthogonal to $P_1=K+E+\sum c_j\tilde{C}_j -N_1$, so is a first step in constructing the negative part of the Zariski decomposition.  If $P_1$ is nef, then $P_1+N_1$ is in fact this decomposition.   Otherwise, one needs to find an $N_2$, as described in Remark $2.2$.

     \begin{theorem} \label{t}  Let $(\tilde{X},\tilde{C}\cup E)$ be the minimal orbifold log resolution of $(X,C=\sum c_iC_i)$, which is not star-shaped.   Suppose that all non-$0$ $c_i$ are $\geq 1/2$.  Then
          \begin{enumerate}
     \item $N\equiv\sum (1-c_i)D_i$ is the negative part of the Zariski decomposition of $K+E+\sum c_i\tilde{C}_i$.
     \item  If $(X,C)$ is log canonical, then the graph of $E$ is
$$
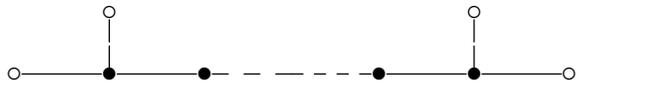
\xymatrix@R=8pt@C=30pt@M=0pt@W=0pt@H=0pt{
&&&\righttag{\circ}{}{6pt}&&&&\righttag{\circ}{}{6pt}&&&\\
&&&\lineto[u]&&&&\lineto[u]&&\\
&&\undertag{\circ}{}{4pt}\lineto[r]&\undertag{\bullet}{}{4pt}\lineto[r]\lineto[u]&\undertag{\bullet}{}{4pt}\dashto[r]&\dashto[r]&\undertag{\bullet}{}{4pt}\lineto[r]&\undertag{\bullet}{}{4pt}\lineto[r]\lineto[u]&\undertag{\circ}{}{4pt}&,\\
&&&&&\\
}
$$  
Each $\circ$ is either a $-1$ curve with an arrow sticking out, corresponding to a $\tilde{C}_i$ with coefficient $c_i=1/2$; or, a $-2$ exceptional curve.

     \end{enumerate}
     \begin{proof} We prove that $P=K+E+\sum c_i\tilde{C}_i-N$ is nef.  By construction, it dots to $0$ with every curve in a string (in particular, $P\cdot N=0$).  Next, consider a curve $F$ not in a string but intersecting $S_1,\cdots, S_t$, with invariants $m_i/q_i, 1\leq i \leq t$.  Assume $F$ has genus $g$ and valency $t'$ in the graph, with $t'>t$ because the star-shaped case is excluded.  Then 
     $$P\cdot F=2g-2+t'-\sum_{i=1}^t(1-c_i)/m_i. $$ If $g>0$, then $P\cdot F>0$. If $t'\geq t+2$, then $P\cdot F\geq 0$.  So consider the case $g=0$ and $t'=t+1$.  Further, $t>1$, else $F$ would be part of a string.  Now examine $$P\cdot F=-1+\sum_{i=1}^t(1-(1-c_i)/m_i).$$  
     The case $m_i=1$ corresponds to a $-1$ string, so $c_i\geq 1/2$ by assumption.  Then for this term, $1-(1-c_i)/1 \geq 1/2$.  If $m_i\geq 2$, then the corresponding term is also $\geq 1/2$.  As $t\geq 2$, one has $P\cdot F \geq 0$.
     
      If $F$ is not in or adjacent to a string, clearly $P\cdot F\geq 0$.  Thus, $P$ is nef.
  
 Assume now that $P=0$.  Returning to $F$ as above intersecting $t$ strings, one must have $g_F=0$.  As noted, each $(1-c_i)/m_i\leq 1/2$, so $P\cdot F=0$ implies $t'< t+2$.  Thus $t'$ must equal $t+1$, and $t\geq 2$ (else $F$ would be part of a string).  One concludes that  $$1=\sum_{i=1}^t(1-(1-c_i)/m_i).$$  
 Since each term in the sum is at least $1/2$, we conclude that $t=2$ and each term equals $1/2$.  Thus, either $m_i=1$ and $c_i=1/2$, or $m_i=2$ and $c_i=0$.  
 
 If $F$ is a curve which is not on a string and does not intersect any, then $P\cdot F=0$ implies that $g_F=0$ and $t'=2$.  Thus the resolution graph is as claimed.
 
 \end{proof}
     \end{theorem}  
     
     \begin{remark}  The $5$ possible configurations of the numbers $-E_j\cdot E_j$ for the curves $E_j=\circ$ are described from left to right via  $1$'s and $2$'s, as $1 1 1 1$; $1 1 1 2$; $1 1 2 2$; $1 2 1 2$; and $1 2 2 2$.
     \end{remark}

  \section{{Maps and log covers of pairs}}

 An orbifold pair $(X,\sum (1/n_i)C_i)$ gives rise on the boundary to an orbifold of the form $\tilde{\Sigma}=(\Sigma, \gamma_1,\cdots,\gamma_r, n_1,\cdots,n_r)$.
 An orbifold cover of $\tilde{\Sigma}$  is a finite map to $\Sigma$ which is unramified except over the $\gamma_i$, over which it is branched with ramification order $\leq n_i$.
 (These facts are recalled in \cite[\S $1$]{nw}). By Definition 2.1 of \cite{nw}, a \emph{log cover of orbifold pairs} $f:(Y,C')\rightarrow (X,C=\sum_{i=1}^r(1/n_i)C_i )$
 is a finite map of normal germs $f:Y\rightarrow X$ so that
\begin{enumerate}
\item $f^{-1}(|C|)=|C'|$, and $f$ is an unramified covering
  off this set.
\item Writing  $f^*C_i=\sum_{j=1}^{s_i}m_{ij} C'_{ij}$, one has $m_{ij}|n_i$, all $j$.
\item $C'=\sum_{i=1}^r \sum_{j=1}^{s_i}(m_{ij}/n_i) C'_{ij}$.

\end{enumerate}
Here $m_{ij}$ is the ramification (or branching) order of $f$ at the generic point of $C'_{ij}$. The conditions $(2)$ and $(3)$ say that $f^*(C)=C'.$  We extend the above definition to
curves $(X,C=\sum c_iC_i)$, where $c_i\in [0,1]\cap \Q$;  if $f:Y\rightarrow X$  has $f^*(C_i)=\sum_{j=1}^{s_i}m_{ij} C'_{ij}$,  define $C'=f^*(C)=\sum c_im_{ij}C'_{ij}$ as before, but only if all coefficients $c_im_{ij}\leq 1$.
\begin{definition} \label{b} A \emph{map of pairs} is a finite morphism 
 $f:(Y,C'=\sum c'_jC'_j)\rightarrow (X,C=\sum c_iC_i)$ for which $f^*(C)=C'$. A \emph{log cover of pairs} is a map of pairs $f:(Y,C')\rightarrow (X,C)$ which is unramified off $f^{-1}(\text{Supp}\  C)$.

\end{definition} 

From a finite map of germs $f:Y\rightarrow X$, there is a natural way to make maps and covers of pairs.
\begin{proposition}\cite[(2.2)]{nw} Suppose $f:Y\rightarrow X$ is finite, and $C_1,\cdots,C_r$ are  irreducible curves on $X$.  Writing $f^*(C_i)=\sum_{j=1}^{s_i}m_{ij} C'_{ij}$ and $n_i=lcm(m_{ij})$, define $C=\sum (1/n_i)C_i$ and $C'=\sum_{i=1}^r \sum_{j=1}^{s_i}m_{ij}/n_i\ C'_{ij}$.  Then
\begin{enumerate}
\item $f:(Y,C')\rightarrow (X,C)$ is a map of pairs.
\item If $f$ is unramified off $f^{-1}(\cup C_i)$, then $f$ is a log cover of pairs.
\end{enumerate}
\end{proposition}

Section $6$ below recalls basic ways to construct log covers of a given pair $(X,C)$.

We need a simple but important definition.  If $C=\sum_{i=1}^r c_iC_i $, with all $c_i>0$, we define the \emph{opposite curve} $\bar{C}$ of $C$ by
$$ \bar{C}:=\sum_{i=1}^r(1-c_i)C_i.$$
 \begin{definition} \label{c} The \emph{volume of a pair} $(X,C)$ is $$Vol(X,C):= -P_{\bar{C}}\cdot P_{\bar{C}}.$$
  \end{definition}
This quantity is computed from a log resolution $(\tilde{X},\tilde{C}\cup E)\rightarrow (X,C)$ via the Zariski decomposition $$K_{\tilde{X}}+E+\sum_{i=1}^r (1-c_i)\tilde{C}_i=P_{\bar{C}}+N_{\bar{C}}.$$
Note that volumes satisfy (using Proposition (3.3)): 
$$(5.3.1)\ C'\leq C\  \text{on}\  X \ \text{implies}\ 
Vol(X,C) \leq Vol(X,C').\ \ \ \ \ \ \ \ \ \ \ \ \ \ \ \ \ \ \ \ \ \ \ $$
$$(5.3.2)\ C\ \text{reduced} \ (\text{all}\ c_i=1)\ \text{implies}\ Vol(X,C)=-P_X \cdot P_X=Vol(X).$$
A basic point is that volumes multiply by degree in a log cover.
     \begin{theorem} \label{g} Suppose $f:(Y,C')\rightarrow (X,C)$ is a log cover of  pairs of degree $d$.  Then $$
      Vol(Y,C')=d\cdot Vol(X,C).$$
      \begin{proof}  We follow the proof in \cite[Theorem 2.8]{w}.  Let $(\tilde{X},\tilde{C}\cup E)\rightarrow (X,C)$ be the minimal log resolution so that every $\tilde{C}_i$ intersects a rational end-curve $E_i$ of $E$, and $E_i$ intersects at most one such curve.  Denote by $\tilde{C}$ the reduced sum of the $\tilde{C_i}$.  Let $f':Y'\rightarrow \tilde{X}$ be the normalization of $\tilde{X}$ in the quotient field of $(Y,0)$.  This finite map is unramified off the inverse image of $E\cup \tilde{C}$, and $Y'$ has only cyclic quotient singularities.  Above the intersection point $P_i=\tilde{C}_i \cap E$ are $s_i$ cyclic quotient singularities (or smooth points), at points $P_{ij}$ on some Weil divisor $\tilde{C}'_{ij}$ on $Y'$.  A local equation defining $\tilde{C}_i$ on $\tilde{X}$ has a zero of multiplicity $m_{ij}$ along $\tilde{C}'_{ij}$.
      
       Let $\tilde{Y}\rightarrow Y'$ be the result of first minimally resolving all the cyclic quotient singularities on $Y'$, and then insuring that the proper transforms $\tilde{C}''_{ij}$ of the $\tilde{C}'_{ij}$ intersect transversally a unique end curve of the reduced exceptional fibre $F$ of the map $\tilde{Y}\rightarrow Y$.  
  
       Letting $\tilde{f}:\tilde{Y}\rightarrow Y'\rightarrow \tilde{X}$ be the induced map, we conclude that $$\tilde{f}^*(\tilde{C}_i)=\sum (m_{ij}\tilde{C}''_{ij}+F_{ij}),$$
       where $F_{ij}$ is a (not necessarily reduced) effective exceptional cycle collapsed to a point on $Y'$.
      
      Since $\tilde{f}$ is proper, generically finite and unramified off the inverse image of the normal crossings divisor $E+\tilde{C}$ , it follows (e.g., \cite[proof of (2.8)]{w}) that $$\tilde{f}^*(K_{\tilde{X}}+E+\tilde{C})=K_{\tilde{Y}}+F+\tilde{C}'',$$ where as before $\tilde{C}''$ denotes the reduced sum of the $\tilde{C}''_{ij}$.  Therefore, 
      $$\tilde{f}^*(K_{\tilde{X}}+E+\sum_{i=1}^r(1-c_i)\tilde{C}_i)=K_{\tilde{Y}}+F+\sum_{i=1}^r\sum_{j=1}^{s_i}(1-c_im_{ij})\tilde{C}''_{ij}-\sum_{i=1}^r\sum_{j=1}^{s_i}c_i F_{ij}.$$  As the left hand side is $\tilde{f}^*(P_{\bar{C}}+N_{\bar{C}})$, we can rewrite as
  $$\tilde{f}^*(P_{\tilde{C}})+ (\tilde{f}^*(N_{\tilde{C}}) +\sum_{i,j} c_iF_{ij})=  K_{\tilde{Y}}+F+\sum_{i,j} (1-c_im_{ij})\tilde{C}''_{ij}.$$
The Zariski decomposition of the right-hand side is $P_{\bar{C'}} +N_{\bar{C'}}$.  We claim that the left-hand side also appears as a Zariski decomposition. First, the purportedly negative part  $\tilde{f}^*(N_{\tilde{C}}) +\sum_{i,j} c_iF_{ij})$ is effective.  Next, $\tilde{f}$ factors through $Y'$, so that $\tilde{f}^*$ of any $\Q$-divisor on $\tilde{X}$ dots to $0$ with all components of each cycle $F_{ij}$.  So, the two summands are orthogonal.   
 Therefore, $\tilde{f}^*(P_{\bar{C}})=P_{\bar{C}'}$.  Self-intersections
 multiply by the degree of $\tilde{f}$ under pull-back, hence the result claimed.

      \end{proof}
       \end{theorem}     

 
  \begin{theorem}\label{h} Suppose $f:(Y,C')\rightarrow (X,C)$ is a map of  pairs of degree $d$.  Then $$
      Vol(Y,C')\geq d\cdot Vol(X,C).$$
      \begin{proof}  We again follow the proof of \cite[Theorem $2.9$]{w}.  Suppose $\tilde{f}:\tilde{Y}\rightarrow \tilde{X}$ is as in the proof of Theorem 5.5 above, but we may also assume that both branch and ramification loci have strong normal crossings divisors.  Let $R\subset \tilde{X}$ denote the reduced sum of those components of the branch locus of $\tilde{f}$ which are not among the curves in $E+\tilde{C}$; let $R'\subset \tilde{Y}$ be the corresponding reduced sum of components of the inverse image of $R$.  Since $\tilde{f}$ is proper and generically finite and unramified off the inverse image of the normal crossings divisor $E+\tilde{C}+R$ , it follows that $$\tilde{f}^*(K_{\tilde{X}}+E+\tilde{C}+R)=K_{\tilde{Y}}+F+\tilde{C}''+R'.$$  Therefore, $$\tilde{f}^*(K_{\tilde{X}}+E+\sum_{i=1}^r(1-c_i)\tilde{C}_i+R)=K_{\tilde{Y}}+F+R'+\sum_{i,j} (1-c_im_{ij})\tilde{C}''_{ij}-\sum_{i,j} c_i F_{ij}.$$
      Note $\tilde{f}^*(R)-R'=R''+Z,$ where $R''$ is effective and supported on $R'$, and $Z$ is an effective exceptional divisor.   Thus 
      $$\tilde{f}^*(K_{\tilde{X}}+E+\sum_{i=1}^r(1-c_i)\tilde{C}_i)+(Z+\sum_{i,j} c_i F_{ij})+R''=K_{\tilde{Y}}+F+\sum_{i,j} (1-c_im_{ij})\tilde{C}''_{ij}.$$
 Consider the following $\Q$-divisors and their Zariski decompositions:$$L=\tilde{f}^*(K_{\tilde{X}}+E+\sum_{i=1}^r(1-c_i)\tilde{C}_i)=\tilde{f}^*(P_{\bar{C}})+\tilde{f}^*(N_{\bar{C}})$$ 
 $$L'=L+(Z+\sum_{i,j} c_i F_{ij})=\tilde{f}^*(P_{\bar{C}})+(\tilde{f}^*(N_{\bar{C}})+Z+\sum_{i,j} c_i F_{ij})$$ 
  $$L''=L'+R''= K_{\tilde{Y}}+F+\sum_{i,j} (1-c_im_{ij})\tilde{C}''_{ij}=P_{\bar{C'}}+N_{\bar{C'}}.$$   
  The Zariski decomposition of $L'$ follows because $Z+\sum_{i,j}c_iF_{ij}$ is effective and is collapsed to points by $\tilde{f}$, hence dots to $0$ with $\tilde{f}^*(P_{\bar{C}})$. 
   
   So,  the $-P\cdot P$ invariant for both $L$ and $L'$ is $$-\tilde{f}^*(P_{\bar{C}})\cdot \tilde{f}^*(P_{\bar{C}})=-(\text{deg}\  \tilde {f})(P_{\bar{C}}\cdot P_{\bar{C}})=d(-P_{\bar{C}}\cdot P_{\bar{C}}),$$   
 while the $-P\cdot P$ invariant for $L''$ is $$-(P_{\bar{C'}} \cdot P_{\bar{C'}}).$$  Replacing a $\Q$-divisor $L$ by an integral multiple $mL$ multiplies the $-P\cdot P$ invariant by $m^2$.  Thus, to prove the desired inequality of the theorem, it suffices to work with $mL'$ and $mL''$ for an integer $m$ making the two $\Q$-divisors integral, thus corresponding to line bundles $\mathcal L'=\mathcal O(mL')$ and $\mathcal L''=\mathcal O(mL'')=\mathcal L'(mR'')$ on $\tilde{Y}$.  Using these line bundles and the result recalled in (2.3.2), it suffices to show that for all $n>>0$,
$$ \text{dim}\ H^0(\tilde{Y}-F,\mathcal L'^n)/H^0(\tilde{Y}, \mathcal L'^n)\leq \text{dim}\ H^0(\tilde{Y}-F,\mathcal L''^n)/H^0(\tilde{Y}, \mathcal L''^n).$$
But this inequality follows from the assertion $$H^0(\tilde{Y}-F,\mathcal L'^n)\cap H^0(\tilde{Y},\mathcal L''^n)=H^0(\tilde{Y},\mathcal L'^n).$$
 This equality is true by local considerations simply because $\mathcal L''$ differs from $\mathcal L'$ by the effective divisor $mR''$ containing no components of $F$; if a section of $\mathcal L'^n$ over $\tilde{Y}-F$ acquired a pole over $F$, it could not extend to a holomorphic section of $\mathcal L''^n$ (which only allows new poles over $R''$).

   \end{proof}  

  \end{theorem}
  

\section{{Log covers and volumes of orbifold pairs}}
In \cite{w}, a germ $(X,0)$ produces an invariant $-P\cdot P$ which on the boundary $\Sigma$ is a strong ``characteristic number'' for $3$-manifolds in the sense of Milnor-Thurston \cite{MT}: multiplicative by degree in unramified coverings, and submultiplicative for arbitrary maps.  This justifies calling it a ``volume.''  Similarly, the boundary of an orbifold pair $(X,\sum(1/n_i)C_i)$ is an orbifold $\tilde{\Sigma}=(\Sigma, \gamma_i,n_i)$, for which the aforementioned $-P_{\bar{C}}\cdot P_{\bar{C}}$ has analogous ``characteristic'' properties.

As a particular case of Definition $5.3$, we have  \begin{definition} $Vol(X,\sum (1/n_i)C_i):=-P_{\bar{C}}\cdot P_{\bar{C}}$, where $\bar{C}=\sum(1-1/n_i)C_i.$
   \end{definition}
   From Theorems $(5.4)$ and $(5.5)$ of the last chapter, we get
   \begin{theorem}\label{k} Suppose $f:(X',\sum (1/n'_j)C'_j)\rightarrow (X,\sum (1/n_i)C_i)$ is a map of degree $d$ of orbifold pairs. Then $$Vol (X',C')\geq d\cdot Vol(X,C),$$ with equality if $f$ is a log cover.
   \end{theorem}
For an orbifold pair $(X,\sum(1/n_i)C_i)$, $(5.3.1)$ implies that if the $n_i$ increase, then the volumes increase (or stay the same). 
  
    There is an analogue of the fundamental group of manifolds for describing orbifold covers of orbifolds (recalled in \cite{nw}, \S$1$).  
\begin{definition}   Consider an orbifold $\tilde{\Sigma}=(\Sigma, \gamma_1,\cdots,\gamma_r,n_1,\cdots,n_r).$ The \emph{orbifold fundamental group} is
$$\pi_1^{orb}(\tilde{\Sigma}):=\pi_1(\Sigma -\bigcup_i\gamma_i)/<\mu_1^{n_1},\cdots,\mu_r^{n_r}>,$$ where $\mu_i$ is represented by the boundary of a
small transverse disc to $\gamma_i$. 
\end{definition}  
One has the usual Galois correspondence between subgroups of finite index of $\pi_1^{orb}(\tilde{\Sigma})$ and finite orbifold covers of $\tilde{\Sigma}$.  There is a \emph{universal orbifold cover} of $\tilde{\Sigma}$ as well 
as a \emph{universal abelian orbifold cover}\  (or UAOC) corresponding to the abelianization of $\pi_1^{orb}$.  Recall the 
\begin{proposition}(\cite[(1.3)]{nw}) Let 
$\tilde{\Sigma}=(\Sigma, \gamma_i, n_i)$ be an orbifold for which $\Sigma$ is a $\Q$HS.  
\begin{enumerate} \item The UAOC\ \  $\tilde{\Sigma}^{ab}$ of $\tilde{\Sigma}$ is a finite orbifold cover.
\item $\tilde{\Sigma}^{ab}=\Sigma^{ab}$ is a manifold, for which all weights $n_j^{ab}=1$.
\item $\pi_1^{orb}(\tilde{\Sigma}^{ab})=\pi_1(\Sigma^{ab})$.
\item The covering group $G=H_1^{orb}(\tilde{\Sigma},\Z)$ sits in an exact sequence $$0\rightarrow \Z/(n_1)\oplus\cdots\oplus\Z/(n_r)\rightarrow G\rightarrow H_1(\Sigma;\Z)\rightarrow 0.$$
\end{enumerate}
\end{proposition}
When the orbifold arises from an orbifold pair $(X,\sum (1/n_i)C_i)$, then log covers give orbifold covers, and one has the 
\begin{theorem}(\cite[(3.1)]{nw}) Suppose $(X,C)$ is an orbifold pair for which the link $\Sigma$ is a $\Q$HS.  Then there exists a finite \emph{universal abelian log cover (or UALC) $(X',C')\rightarrow (X,C)$}, with $C'$ reduced, which induces the UAOC on the boundary.
\end{theorem}
Since $C'$ is reduced in the last result, $Vol(X',C')=Vol(X')$, and we get
\begin{corollary}  Suppose $(X,\sum (1/n_i)C_i)$ is an orbifold pair for which the link $\Sigma$ is a $\Q$HS, and let $f:(X',C')\rightarrow (X,C)$ denote the UALC.  Then $$Vol(X')=Vol(X,C)\cdot |H_1(\Sigma,\Z)|\cdot \Pi n_i.$$
\end{corollary} 
We are now ready to characterize topologically the orbifold pairs of volume $0$, analogous to the aforementioned theorem that $(X,0)$ is log canonical if and only if it has volume $0$ if and only if its local fundamental group is finite or solvable.  The complete list of these pairs is given in the next Section.
\begin{theorem} \label{l}  An orbifold pair $(X,\sum (1/n_i)C_i)$ has volume $0$ if and only if the corresponding orbifold fundamental group is finite or solvable.
\begin{proof}  Suppose $(X,\sum (1/n_i)C_i)$ has $\Q$HS link.  Then by the Corollary it has volume $0$ iff the UALC $(X'',C'')$ has volume $0$, which is the same as $Vol(X'',0)=0$.  Further, since $\pi_1^{orb}(\tilde{\Sigma''})=\pi_1(\Sigma'')$ is a finite index normal subgroup of $\pi_1^{orb}(\tilde{\Sigma})$, the finiteness
 or solvability of one of these groups is equivalent to that property for the other.  But $\pi_1(\Sigma'')$ is finite or solvable if and only if the volume of $(X'',0)$ is $0$.  This proves the result in these cases.  
 
  If $(X,C)$ has volume $0$, then 
$(X,\sum (1-1/n_i)C_i)$ is log canonical, hence ($3.4.1$) $(X,0)$ is log terminal, hence a quotient singularity; one is in the $\Q$HS case, so the result about the orbifold fundamental groups follows.

If $\pi_1^{orb}(\tilde{\Sigma})$ is finite or solvable, the same is true for its quotient $\pi_1(\Sigma)$, so $(X,0)$ is log canonical.
The only log canonical $(X,0)$ which don't have $\Q$HS links are simple elliptic and cusp singularities (which have infinite solvable local fundamental groups).  Rather than handle these two special cases, we prove a more general result below. It asserts that every orbifold pair $(X,C)$ has a log cover $(X',C')$, with $C'$ reduced.  Then $\pi_1^{orb}(\tilde{\Sigma'})=\pi_1(\Sigma')$ is a finite index subgroup of $\pi_1^{orb}(\tilde{\Sigma})$, hence is itself finite or solvable.  Therefore, $(X',C')$ has volume $0$, and thus so does $(X,C)$.
\end{proof}
\end{theorem}

As in \cite{nw}, adjoining $n^{th}$ roots and normalizing gives cyclic log covers. 
\begin{proposition}(\cite[(2.4)]{nw})  Let $C\subset X$ be a reduced Weil divisor, $h$ a regular function with zero divisor $(h)=C+nD$, for some effective $D$. Then adjoining an $n^{th}$ root of $h$ to the function field of $(X,0)$ and normalizing gives a log cover $(X',C')\rightarrow (X,(1/n)C)$, where 
  $C'$ is reduced.  
\end{proposition}
 Imitating the proof from \cite{nw} of Theorem $6.5$ above, one can drop the $\Q$HS restriction and prove the 

\begin{theorem}  For any orbifold pair $(X,C=\sum (1/n_i)C_i)$, there exists a log cover $f:(Y,C')\rightarrow (X,C)$ with $C'$ reduced.  In particular, $$Vol(X,C)=Vol(Y)/\text{deg}(f).$$
\begin{proof}  We rely on the following (Cf. \cite[(2.2)]{o}):
\begin{lemma}  For any $g$ in the divisor class group $D=\text{Cl}(X,0)$ and $n\in \N$, there exist $t\in \text{Tor}(D)$ and $u\in D$ so that $g=t+nu.$
\begin{proof} For any good resolution $(\tilde{X},E)\rightarrow (X,0)$, writing $\E=\bigoplus \Z\cdot E_i$ one has the well-known exact sequence of Mumford (\cite[p. 16]{M}):
\begin{center} $H^1(\tilde{X},\Z)\rightarrow H^1(\tilde{X},\mathcal O_{\tilde{X}})\xrightarrow{\alpha} D \xrightarrow{\beta} \E^*/\E \rightarrow 0.$
\end{center}  As $\beta(g)$ is torsion,  $\beta(mg)=0$ for some $m$,  hence one can write $mg=\alpha(v)$. As  $v\in H^1(\tilde{X},\mathcal O_{\tilde{X}})$, a complex vector space, one can define $t=g-\alpha((1/m)v)\in D$, a torsion element.  Then $g=t+n\alpha((1/mn)v).$
\end{proof}
\end{lemma}
Returning to the Theorem, write the class of $C_i$ as $[C_i]=t_i+n_i[D_i]$, for torsion elements $t_i$ and classes $[D_i]$.  Use the $t_i$ to successively form unramified coverings of $(X,0)$, giving an unramified map $g:(X',0)\rightarrow (X,0)$ for which each $g^*(t_i)$ is trivial in $\text{Cl}(X',0)$.  Thus $g^*[C_i]-n_ig^*[D_i]$ is trivial, hence $g^*(C_i)-n_ig^*(D_i)$ is a principal divisor.  As in the proof of Theorem $3.1$ of \cite{nw}, there is a regular function $h_i$ on $X'$ defining an effective divisor of the form $g^*(C_i)+n_iF_i$, for some curve $F_i$.  Then, successively taking $n_i^{th}$ roots of the functions $h_i$ and normalizing gives the desired result.  
\end{proof}
\end{theorem}  
\begin{remark}  The decomposition in Lemma $6.10$ is not unique when $(X,0)$ does not have $\Q$HS link, hence the log cover in Theorem $6.7$ need not be unique (nor a Galois cover).
\end{remark}

  \section{{Orbifold pairs with volume 0}}
   
 An orbifold pair $(X,\sum (1/n_i)C_i)$ has volume $0$  exactly when $(X,\sum (1-1/n_i)C_i)$ is log canonical.  Theorem $4.6(2)$ gives the only examples for which the minimal log resolution is not star-shaped.  (A coefficient $c_i=1/2$ means orbifold weight $2$.)  Lemma $4.1$ covers rational chains (they are all log canonical, with arbitrary weights on the ends).  It  only remains to consider star-shaped graphs with $t\geq 3$ branches.
 
 Assume the determinants of the branches are $m_1,m_2,\cdots,m_t$, with weights $n_1,\cdots,n_t$ (recall we set $n_j=1$ if there is no curve $\tilde{C_j}$ on the minimal log resolution.)  Then $c_i=1-1/n_i$, and the inequality $\chi_C\leq 0$ of Proposition $4.2$ translates easily to $$\sum_{i=1}^t(1-(1/n_im_i))\leq 2.$$ 
  Note $m_i=1$ implies $n_i\geq 2$ , by minimality of the resolution; so $n_im_i\geq 2$, hence summands are at least $1/2$, whence $t\leq 4$.  One easily concludes:

  \begin{proposition}  The volume $0$ orbifold pairs with $t=4$ have $n_im_i=2$ for all $i$, with graph
$$
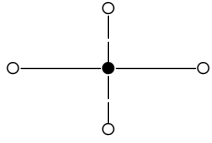
\xymatrix@R=8pt@C=30pt@M=0pt@W=0pt@H=0pt{
&&&\righttag{\circ}{}{6pt}&&&&&&&\\
&&&\lineto[u]&&&&&\\
&&\undertag{\circ}{}{4pt}\lineto[r]&\undertag{\bullet}{}{4pt}\lineto[r]\lineto[u]&\undertag{\circ}{}{4pt}&\\
&&&\lineto[u]&&&&\\
&&&\lineto[u]\righttag{\circ}{}{6pt}&&&\\
&&&&&\\
}
$$  
Each $\circ$ is either a $-1$ curve with an arrow sticking out, corresponding to a $\tilde{C}_i$ with weight $n_i=2$; or, an unweighted $-2$ exceptional curve.
  \end{proposition}
   One is left with the case $t=3$, for which the above inequality becomes $$\sum_{i=1}^3 1/n_im_i\geq 1.$$ Therefore,
  \begin{proposition} The volume $0$ orbifold pairs with $t=3$ are those whose ordered set of triples $\{n_1m_1\leq n_2m_2\leq n_3m_3\}$ equals one of the following: 
  \begin{itemize}
  \item  $\{2\leq 2\leq k\}, k\geq 2$
  \item  $\{2\leq 3\leq k\}, 3\leq k \leq 6$
  \item  $\{2\leq 4\leq 4\}$
  \item $\{3\leq 3\leq 3\}$
  \end{itemize}
  \end{proposition}
 In enumerating all the examples of the $n_i$ and $m_i$ from the list above, recall that at least one $n_i\geq 2$; an $m_i>2$ allows for several chains $m_i/q_i$; self-intersection of the central curve is arbitrary as long as the graph is negative-definite.
 \begin{example}
 For the ordered triple $\{2,3,4\}$, we specify all possible ordered sets of weights $(n_1,n_2,n_3)$ (from which the $m_i$ are determined).  These are: $(1,1,4)$, $(1,1,2)$,
 $(1,3,1)$, $(2,1,1)$, $(1,3,4)$, $(1,3,2)$,
 $(2,1,2)$, $(2,1,4)$, $(2,3,1)$, $(2,3,2)$, $(2,3,4)$.  (Alternatively, one could list possible triples of $m_i$'s.)
 \end{example}
 We summarize the previous results as follows:
 \begin{theorem} \label{m} Suppose $(X,\sum (1/n_i)C_i)$ is an orbifold pair of volume $0$.  Then the graph of the minimal log resolution is one of the following:
 \begin{enumerate}
 \item a string of rational curves  in Lemma 4.1, with arbitrary weights on ends
 \item the graph in Corollary 4.3, with arbitrary weight on the end
  \item the graph in Theorem 4.6(2)
 \item the star-shaped graph of valency four in Proposition 7.1
 \item any star-shaped graph of valency three whose multiplicities and weights are given in Proposition 7.2.
 \end{enumerate}
 \end{theorem}

 Koll\'{a}r (\cite[p. 125]{k}) lists all log canonical pairs $(X,\sum c_iC_i)$ for which all $c_i$ are $0$ or $\geq 1/2$.  If all $c_i$ are of the form $1-1/n_i$, we say the pair comes from a volume $0$ orbifold; these have all been enumerated in Theorem 7.4.  The few examples with a $c_i=1$ ((4.1) and (4.3)), which correspond to orbifold cases for which $\tilde{C}_i$ can have arbitrary weight, can be set aside.
 \begin{lemma}  Let $(X,\sum c_iC_i)$ be a log canonical pair for which all $c_i$ are $0$ or $\geq 1/2$.  Write $c_i=1-1/\bar{n_i}$, where $\bar{n_i}$ equals $1$ or is $\geq 2$. Let $n_i=\lfloor \bar{n_i} \rfloor$.  Then $(X,\sum (1/n_i) C_i)$ is an orbifold pair of volume $0$, with the same graph but smaller weights.
 \begin{proof}   $(X,\sum (1-1/n_i) C_i)$ is less than or equal to  $(X,\sum c_iC_i)$, hence is itself a log canonical pair.
 \end{proof}
 \end{lemma} 
  \begin{lemma} Let $(X,\sum c_iC_i)$ be a log canonical pair for which all $c_i$ are $0$ or $\geq 1/2$.  Unless the minimal log resolution is star-shaped with $3$ branches, it is a volume $0$ orbifold pair.
  \begin{proof} The preceding lemma shows any such pair with possibly lowered weights is an orbifold pair.  But Theorem 7.4 shows that examining all cases except for the star-shaped one with $t=3$, the allowed $c_i$ for the graphs are $0$ or $1/2$.
  \end{proof}
 \end{lemma}
 Lowering the weights $n_i$ on a star-shaped $t=3$ case but leaving the $m_i$'s alone, one uses Proposition $7.2$ and Lemma $7.5$ to conclude:
 \begin{proposition}  Suppose $(X,\sum c_iC_i)$ is a log canonical pair with all non-$0$ $c_i\geq 1/2$, whose minimal orbifold resolution is star-shaped with three branches.  The possibilities for the determinants $(m_1,m_2,m_3)$ of the branches are:
 \begin{enumerate}
 \item $\{1,1,1\}$
 \item $\{1,1,m\}, m\geq 2$
 \item $\{1,2,m\}, m\geq 2$
 \item $\{1,3,m\}, m\in \{3,4,5,6\}$
 \item $\{1,4,4\}$
 \item $\{2,2,m\}, m\geq 2$
 \item $\{2,3,3\}$
 \end{enumerate}
 \end{proposition}
 \begin{example}  Consider $(\C^2, c_1L_1+c_2L_2+c_3L_3)$, where the $L_i$ are lines through the origin.  The minimal log resolution is star-shaped with a central $-4$, and the three branches are $-1's$ with weights the $c_i$.  If all coefficients $c_i$ are $\geq 1/2$, it is log canonical iff $c_1+c_2+c_3\leq 2$.  The orbifold examples are clearly sparsely distributed among these.
 \end{example}
 
 To work out each of the other examples, one should separate cases based on whether an $m_i>1$ corresponds to a $c_i=0$ or $c_i\geq 1/2$.

 Since $\overline{1/2C}=1/2C$ for a reduced curve, one can ask for which
 plane curves $C$ is $Vol(\C^2,(1/2)C)=0$.   
 \begin{proposition}  Let $C=\{f(x,y)=0\} $ be reduced.  Then $(\C^2,(1/2)C)$ has volume $0$ if and only if $f$ is  one of the following (up to equisingularity):
 \begin{enumerate}
 \item $y$
 \item $xy$
 \item $y^2-x^{2n}$, \ $n\geq 2$
 \item $y^2-x^{2n+1}$,\ $n\geq 1$
 \item $x(y^2-x^{2n})$,\ $n\geq 2$
 \item $x(y^2-x^{2n+1})$, $n\geq 2$
 \item $y^3+x^4$
 \item $y(y^2+x^3)$
 \item $y^3+x^5$
 \item $y^3+x^6$
 \item $y^4+x^4$
 \item $(y+x^2)(y^2+x^{2n+4})$
 \item $(y+x^2)(y^2+x^{2n+3})$,
 \item $x(y+x)(y^2+x^{2n})$
 \item $x(y+x)(y^2+x^{2n+1})$
 \end{enumerate}
 \begin{proof}  By Proposition 6.8, the Volume of the double cover $\{z^2=f(x,y)\}$ is twice the Volume of the pair $(\C^2,(1/2)C)$, hence one is log canonical if and only if the other is.  The graphs of log canonical singularities have been well-known since \cite{K}, and they are all rational or minimally elliptic.  So one concludes that the log canonical double points in $\C^3$ are either:  rational double points; two types of simple elliptic singularities (which have a modulus); two types of cusp singularities.  Equations for the cusps (the last four above) can be found on page 1290-91 of \cite{laufert}; the  dual graph type is there called \emph{No}.  

 \end{proof}
 \end{proposition}

\begin{section}{{The set of volumes of RDP orbifold pairs}}

 As remarked in the Introduction, F. Ganter proved in \cite{G} that the set of volumes of Gorenstein NSS's satisfies the DCC, with minimum non-$0$ value $1/42$ (achieved for the Brieskorn singularity $V(2,3,7)$).  However, there is no simple generalization for Gorenstein orbifold pairs.   Example 8.8 below gives a family $(X_m,C_m)$ of Gorenstein (even minimally elliptic) pairs for which the volumes $Vol(X_m,C_m)$ are strictly decreasing (here, the volumes $Vol(X_m)$ are strictly increasing).  It therefore makes sense to formulate a version for which the $X$ are Gorenstein log terminal, i.e. RDP's.
 
\begin{DVC}  The set of volumes of RDP orbifold pairs $(X,C=\sum (1/n_i)C_i)$ satisfies the DCC, with minimum non-$0$ value $1/3528=1/(2\cdot 42^2)$.
\end{DVC}

\begin{example} The orbifold pair $$(X,C)=(\C^2,(1/2)\{x=0\}+(1/3)\{y=0\}+(1/7)\{x+y=0\})$$ has minimal log resolution graph
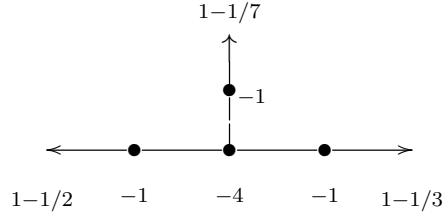
$$
\xymatrix@R=8pt@C=30pt@M=0pt@W=0pt@H=0pt{
&&&&\overtag{}{1-1/7}{6pt}&&&&&\\
&&&&\uparrow&&\\
&&&&\righttag{\bullet}{-1}{6pt}\lineto[u]&&&&\\
&&&&\lineto[u]&&&\\
&&\undertag{~}{1-1/2}{12pt}&\ar[l]\undertag{\bullet}{-1}{6pt}\lineto[r]&\undertag{\bullet}{-4}{6pt}\lineto[r]\lineto[u]&\undertag{\bullet}{-1}{6pt}\lineto[r]&\undertag{\ar}{1-1/3}{12pt}&\\
&&&&&&\\
&&&&&\\
}
$$ 
where the arrows correspond to the curves in $C$.
 The UALC is given by adjoining $u,v,w$ satisfying $u^2=x,v^3=y,w^7=x+y$; this is the hypersurface $\{u^2+v^3=w^7\}$, of volume $1/42$.  Therefore, $Vol(X,C)=1/(2\cdot 3\cdot7)\cdot (1/42).$
Since $\{\pm \text{Id}\}$ acts freely on $(X,C)$, the quotient $(Y,C')$ is an $A_1$ singularity with weights along three lines.  The minimal log resolution is the same as the one above, except that the self-intersection of the central curve is $-5$. One has $Vol(Y,C')=(1/2)Vol(X,C)=(2\cdot 42\cdot 42)^{-1}=1/3528.$

\end{example}
We will verify the Conjecture in the star-shaped case (as done in \cite[(3.5)]{w} for volumes of singularities, before \cite{G} proved the general case).
\begin{theorem} \label{n} The set of volumes of RDP orbifold pairs with star-shaped minimal log resolution satisfies the DCC, with minimum non-$0$ volume $1/3528.$
\begin{proof}
Here the log resolution graph has some number $s$ of $-1$-curves attached to a central vertex, and blowing those down gives an RDP, with $r\leq 3$ branches emerging from the center.  The new central vertex is either a $-2$ (in which case the graph now is a familiar ADE) or a $-1$ (which can happen only for $r\leq 2$).  The RDP consists either of $r$ chains of determinant $m_1,\cdots,m_r$ connected to the central vertex, with $r$ is $1, 2$, or $3$; or ($r=0$) one simply has $s\geq 3$ $-1$-curves emerging from a central $-(s+1)$ or $-(s+2)$.  To exclude volume $0$ examples, one has to suppose that $r+s\geq 3$.

 Write the orbifold pair as $(X,C=\sum_{i=1}^r(1/n_i)C_i +\sum_{j=1}^s (1/n_j')C_j')$, where $C_i$ is at the end of the chain corresponding to $n_i$, and each $C'_j$ is at an end of a $-1$ curve emerging from the central vertex.  Note $n_i\geq 1$, $m_i\geq 2$ and $n_j'\geq 2$.   One has the negative quantity $\chi=-2+\sum_{i=1}^r(1-1/m_i)$, while  positive volume of the pair requires the positivity of $\chi_{\bar{C}}=\chi+\sum_{i=1}^r(1-1/n_i)/m_i+\sum_{j=1}^s(1-1/n'_j)$.   We rewrite it as 
\begin{equation}\label{*}\tag{$*$}
\chi_{\bar{C}}=-2+\sum_{i=1}^r(1-1/n_im_i)+\sum_{j=1}^s(1-1/n'_j).
\end{equation}
\begin{lemma} $\chi_{\bar{C}}>0$ implies $\chi_{\bar{C}}\geq 1/42$, achieved when $r+s=3$, all $n_i=1$, and $\{m_i, n_j'\}=\{2,3,7\}$.
\begin{proof} The expression $-2+\sum_{i=1}^u(1-1/k_i)$, where all $k_i\geq 2$, is clearly negative if $u\leq 2$ and $\geq 1/2$ if $u\geq 5$.  For $u=4$, the minimum positive value is $1/6$, occurring for the $k_i's$ equal $2,2,2,3$.  For $u=3$, the expression is equal to $1-1/k_1-1/k_2-1/k_3$, and it is easy to check the smallest positive value is $1/42$, for the $k_i's$ equal to $2,3,7$.  This result then implies the assertion of the Lemma.
\end{proof}
\end{lemma}
\begin{lemma} The minimum positive value of $Vol(X,C)$ is $1/3528$.
\begin{proof} The volume equals $\chi_{\bar{C}}^2/e$, where $e=|\chi/t|$ is computed from the entries in the chains and the self-intersection of the central vertex.  By Lemma 8.3, the numerator is $\geq (1/42)^2$.  The denominator is $\leq |\chi|\leq 2$,  so the whole expression is $\geq (1/42)^2/2=1/3528$.  Example 8.1 shows that this case does occur.
\end{proof}
\end{lemma}
 
\begin{lemma}  Given an infinite sequence of graphs whose volumes are strictly decreasing, there is an infinite subsequence with fixed $s$ and with fixed $r$.
\begin{proof}  There are only $4$ values for $r$, so at least one of them occurs infinitely often.  Since $e\leq 2$, $Vol(X,C)\geq\chi_{\bar{C}}/2\geq (-2+s/2)/2$, so $s$ must be bounded in the sequence.  Therefore, at least one value occurs infinitely often, and we take such a subsequence.
\end{proof}
\end{lemma}
To conclude the DCC among a set of volumes, use the likely well-known
\begin{proposition}  Suppose $\N^k$ is given a partial order using the ``reciprocal'' weight function $W(n_1,n_2,\cdots,n_k)=\sum_{i=1}^k 1/n_i$.  Then $\N^k$ satisfies the ACC (ascending chain condition): every subset $Z\subset \N^k$ contains an element of maximal weight.
\begin{proof} The case $k=1$ being obvious, proceed by induction on $k$.  For $1\leq i \leq k$, denote by $\pi_i:\N^k\rightarrow \N^{k-1}$ projection off the $i^{th}$ component.   By induction, $\pi_i(Z)$ contains an element of maximal weight; let $\alpha^i=(\alpha^i_1,\alpha^i_2,\cdots,\alpha^i_k)\in Z$ be an element which maps to it.

The assertion of the proposition is proved if there are only finitely many elements of $Z$ whose weight is greater than or equal to all $W(\alpha^i)$.   Suppose  $\beta\in Z$ has $W(\beta)\geq W(\alpha^i)$ for all $i$, and write $\beta=(\beta_1,\cdots,\beta_k)$.  We claim that $\beta_i\leq \alpha^i_i$, for all $i$; this would imply the finiteness assertion.   Now, $$W(\beta)=1/\beta_i+\sum_{j\neq i} 1/\beta_j \geq W(\alpha_i)=1/\alpha^i_i+\sum_{j\neq i}(1/\alpha^i_j).$$ By the maximality assumption for $\alpha^i$, the sum for $j\neq i$ of $1/\beta_j$ is less than or equal to the corresponding sum of $1/\alpha^i_j$, hence $1/\beta_i\geq 1/\alpha^i_i$, as claimed.
\end{proof}
\end{proposition}

\begin{lemma}  For fixed $r$ and $s$, the set of positive expressions $\chi_{\bar{C}}$ satisfies the DCC.
\begin{proof}  By $(*)$, the $\chi_{\bar{C}}$ is $r+s-2$ minus a sum of $r+s$ reciprocals of integers, so apply the Proposition.
\end{proof}
\end{lemma}
If volumes $\chi_{\bar{C}}^2/e$ are strictly decreasing, to use Lemma 8.7 requires understanding the $e$'s in the sequence.
We handle separately the cases of various $r$.  

For $r=0$, one has $\chi=-2$ and $e=1$ or $2$; one of these cases occurs infinitely often, so Lemma 8.7 can be used.

For $r=3$, we break up into the cases of the $n_i$'s.  For $(2,3,3), (2,3,4), (2,3,5)$, there is only one $e$, so the result follows from Lemma 8.7. For $(2,2,n)$, one has $e=1/n$.  But since $\chi_{\bar{C}}\geq 1/42$, in an infinite decreasing family of this type we have $Vol(X,C)\geq n\cdot (1/42)^2$, whence $n$ must be bounded above.  So, restrict to a subsequence with a fixed $n$, hence fixed $e$, and then apply Lemma 8.7 to the resulting expression.

 In case $r=2$, one has $\chi=-(1/n_1+1/n_2)$ (with $n_i\geq 2$ by minimality) so $1/e=|t|/(1/n_1+1/n_2)$. In this case $Vol(X,C)\geq |t|\cdot 1/3528$, whence $|t|$ is bounded in any decreasing family; restricting to a subsequence we can assume it is fixed, so that one is examining the volumes $\chi_{\bar{C}}^2/(1/n_1+1/n_2)$.  In the infinite sequence, the expression $1/n_1+1/n_2$ is bounded above and below, hence there is an infinite subsequence which is strictly monotone or constant.  Pass to this subsequence in the family of decreasing volumes.  By Proposition 8.6, the subsequence cannot be strictly increasing.  If the sequence becomes constant, then one has $e$ fixed in the family; but Lemma 8.7 then implies a contradiction to an infinite sequence of decreasing volumes.  Finally, if the subsequence is strictly decreasing, then the numerator $\chi_{\bar{C}}^2$ must be strictly decreasing in order for the quotient to  be strictly decreasing. Now apply Lemma 8.7.
 
 The same argument applies in the case $r=1$; $|\chi|$ in the denominator of the volume has the form $1+1/n$ or $1+1/(2n+1)$, so in a subsequence its value either increases or decreases or stays the same.
\end{proof}
\end{theorem}

\begin{example}  Let $X_m$ denote (for $m\geq7$) the quasi-homogeneous minimally elliptic singularity whose minimal good resolution graph has a $-1$ rational central curve surrounded by $3$ rational curves, of self-intersection respectively $-2$, $-3$, $-m$,  with $m\geq 7$.  Then $\chi_m=e_m=1/6-1/m$, hence $Vol(X_m)=1/6-1/m$ strictly increases as $m$ increases. Next let $C_m$ be a curve of weight $2$ on $X_m$ for which $\tilde{C_m}$ is a smooth curve on the resolution passing transversally through the central curve at a fourth point.  The orbifold log resolution (with an unmarked bullet a $-2$) is
$$
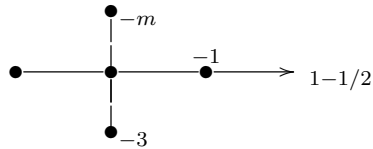
\xymatrix@R=8pt@C=30pt@M=0pt@W=0pt@H=0pt{
&&&&\righttag{\bullet}{-m}{6pt}&&&&\\
&&&&\lineto[u]&&&\\
&&&\overtag{\bullet}{}{6pt}\lineto[r]&\undertag{\bullet}{}{6pt}\lineto[r]\lineto[u]\lineto[d]&\overtag{\bullet}{-1}{6pt}\lineto[r]&\arrow\righttag{}{1-1/2}{6pt}\\
&&&&\lineto[u]&&\\
&&&&\righttag{\bullet}{-3}{6pt}\lineto[u]&&\\
&&&&&\\
}
$$ 
Therefore $\chi_{\bar{C_m}}=\chi_m+(1-1/2)$, and so $Vol(X_m,C_m)=(\chi_m+1/2)^2/(\chi_m)$, which can be rewritten as \begin{center} $8/3+2(4m+3)/m(m-6).$
\end{center}  We see that the volumes of  the pairs $(X_m,C_m)$ strictly decrease as $m$ increases.
\end{example}
\end{section}

\end{document}